\documentclass[final]{amsproc}

\usepackage{showkeys}
\usepackage[english]{babel}
\usepackage{dblaccnt}
\usepackage{accents}
\usepackage{graphicx}
\usepackage{color}
\usepackage{float}
\usepackage{amscd}
\usepackage[mathscr]{euscript}
\usepackage{lipsum}
\usepackage{epstopdf}
\usepackage{algorithm}
\usepackage{algpseudocode}

\usepackage{graphicx,subcaption}
\graphicspath{{images/}}

\newtheorem{theorem}{Theorem}[section]
\newtheorem{lemma}[theorem]{Lemma}

\theoremstyle{definition}

\theoremstyle{remark}
\newtheorem{remark}[theorem]{Remark}

\newcommand{\R}{\mathbb{R}}

\renewcommand{\S}{\mathbb{S}}

\renewcommand{\phi}{\varphi}
\newcommand{\dd}{\mathrm{d}}
\renewcommand{\O}{\mathcal{O}}

\renewcommand{\sc}{\mathrm{sc}}

\renewcommand{\Re}{\mathrm{Re}\,}
\renewcommand{\Im}{\mathrm{Im}\,}

\numberwithin{equation}{section}

\begin{document}

\title[Sampling  method  and   deep learning for  inverse scattering]{Sampling type method  combined with deep learning  for   inverse  scattering  with one incident wave }

\author{Thu Le}
\address{Department of Mathematics, Kansas State University, Manhattan, Kansas 66506}
\email{thule@ksu.edu}
\thanks{This work  was  supported in part by NSF Grant DMS-1812693.}

\author{Dinh-Liem Nguyen}
\address{Department of Mathematics, Kansas State University, Manhattan, Kansas 66506}
\email{dlnguyen@ksu.edu}

\author{Vu Nguyen}
\address{Department of Computer Science, University of Science, Vietnam National  University, Ho Chi Minh City, Vietnam}
\email{vunguyenthai73@gmail.com}

\author{Trung Truong}
\address{Department of Mathematics, Kansas State University, Manhattan, Kansas 66506}
\email{trungt@ksu.edu}




\begin{abstract}
We consider the inverse   problem of 
determining the geometry of penetrable  objects from scattering data generated by  one incident wave at a fixed frequency. We first study an orthogonality sampling type method 
which   is fast, simple to implement, and robust against noise in the data. 
This sampling method has a new imaging functional that is applicable to data measured in  near field or far field regions.
The resolution analysis of the imaging functional is analyzed where the explicit decay rate of the functional is established. 
A connection with the orthogonality sampling method by Potthast is also studied.
The sampling method is then combined with a deep neural network to solve the inverse scattering problem.  
This combined  method can be understood  as a network using
the  image computed by the  sampling method   for the first layer and followed
by the U-net architecture for the rest of the layers.  The fast computation and the  knowledge from the results of the sampling  method help
speed up the training of the network.  The combination leads 
to a significant improvement in the reconstruction results initially obtained by the sampling method.  
The combined method is also  able  to invert some limited aperture experimental data  without any additional transfer training. 
\end{abstract}

\maketitle

%


\section{Introduction}

We consider  an inhomogeneous medium that 
fills up a bounded Lipschitz domain $D \subset \mathbb{R}^{n}$ ($n = 2$ or 3).  Suppose that 
the medium is 
characterized by  bounded function $\eta(y)$ satisfying $\eta = 0$
in $\R^n \setminus \overline{D}$. Consider the incident plane wave
$$
u_{\mathrm{in}}(x) = e^{ikx\cdot d}, \quad  x \in \R^n, \quad d \in  \mathbb{S}^{n-1}: = \{x \in \R^n: |x| = 1\},
$$
where $k>0$ is the wave number and $d$ is the direction vector of propagation.   
The scattering of  $u_{\mathrm{in}}$ by the  inhomogeneous medium is described by the following model problem
\begin{align}
\label{Helm}
& \Delta u + k^{2}(1 + \eta(x)) u =0,\quad  x \in \mathbb{R}^{n},  \\
& u = u_\sc + u_\mathrm{in}, \\
 \label{radiation}
& \lim_{r\rightarrow \infty }r^{\frac{n-1}{2}}\left( \frac{\partial u_{\mathrm{sc}}}{
\partial r}-iku_{\mathrm{sc}}\right) =0,\quad r=|x|,
\end{align}
where $u$ is the total field, $u_\sc$ is the scattered field, and the Sommerfeld radiation condition~\eqref{radiation} holds
uniformly for all  directions $x/|x| \in \S^{n-1}$. 
If $\mathbb{R}^n \setminus \overline{D}$ is connected and $\Re(\eta) \geq 0$,  this scattering problem is known to have a unique weak solution $u_\sc \in H^1_{\mathrm{loc}}(\mathbb{R}^n)$, see~\cite{Colto2013}. 
 
\textbf{Inverse problem.} Let $\Omega$ be a bounded Lipschitz domain such that $\overline{D} \subset \Omega$
and denote by $\nu(x)$ the outward normal unit vector to $\partial \Omega$ at $x$.
Given  $u_\sc$  and $ \partial u_\sc/ \partial \nu$  on $\partial \Omega$ determine  $D$.

In this paper we  develop a sampling type method and combine this method with deep learning to solve the inverse problem above. 
There has been a large body of literature  on both theoretical and numerical studies on this inverse problem and its variants, see~\cite{Colto2000c, Colto2013} and references therein. 
To numerically solve the inverse problem sampling methods are known to be fast, non-iterative  and  do not require advanced a priori information about the 
unknown target. More details about sampling methods  can be found in~\cite{Potth2006, Kirsc2008,   Cakon2011} and references therein.

The sampling method studied in this paper is inspired by  the orthogonality sampling method (OSM)  by Potthast~\cite{Potth2010}. The OSM  is also studied in other works under the name of direct sampling method.  The OSM is  attractive and promising thanks to its computational efficiency. For instance, the OSM is a regularization free method that is very robust against noise in the data,  and its implementation   only involves an evaluation of vector  products. However, the theoretical analysis of the OSM is far less developed compared with that of  the linear sampling method and the factorization method. Also most of the published results deal with the case of far-field measurements, see e.g.,~\cite{Potth2010,Gries2011, Ito2012, Ahn2020, Ito2013, Nguye2019, Harri2020, Le2022}. There have been only a few results on the OSM concerning the case of near-field measurements. The near-field OSM studied in~\cite{Akinc2016} is  applicable to the 2D case with circular measurement boundaries. The 3D  case was studied in~\cite{Kang2021} under the small volume hypothesis of well-separated inhomogeneities. The sampling method proposed in this paper can be applied to near-field data (or far-field data) and it is  not limited to those conditions. However, it requires Cauchy data instead of  only scattered field data.  
Using the Lippmann-Schwinger equation and the Helmholtz integral representation we show that the  imaging functional is equal to a functional that involves Bessel functions $J_0(k|y-z|)$ or $j_0(k|y-z|)$  where $z$ is a sampling point and $y$ is inside the unknown scatterer. These Bessel functions are crucial for  justifying the behavior of the imaging functional.

It is known that in the case of scattering data generated by one incident wave the OSM can only  provide reasonable reconstructions for small scattering  objects~\cite{Potth2010}. It is 
also the case of  other reconstruction methods. This motivates us to combine our orthogonality sampling type method with deep  learning to obtain better reconstruction results. 
The orthogonality sampling type method is combined with deep learning in the following way. 
This combined  method can be understood  as a network that uses
the  image computed by the  sampling method   for the first layer,  followed
by the U-net structure for the remaining layers. The combination  of deep learning
to physics-based reconstruction methods to solve inverse problems has been studied, 
see, e.g.,~\cite{Jin2017, Hamil2018, Li2019, Antho2019, Khoo2019, Xiao2020, Chen2020}. In particular, we refer to the review paper~\cite{Chen2020} 
and references therein for recent advances in using deep learning for solving inverse scattering problems. 
To our knowledge (see also~\cite{Chen2020}) most of the studies on deep learning for inverse scattering consider  data generated by multiple 
incident fields while the the data for inverse scattering problem considered in this paper
 is associated with only one incident field at a fixed wave number. Our inverse problem has a minimal amount of data. 
 In addition,   the imaging functional for the proposed sampling method is 
a new functional and it has an advantage over the imaging functional of the OSM that it can be applied to  near field data or far field data.  
Furthermore, to our knowledge the combination of this sampling method (or orthogonality sampling  methods) and deep learning 
for solving inverse scattering problems has not been studied before.

In this paper we exploit a deep neural network of U-Net-Xception style~\cite{Choll2016}. We train this network using simulated data sets of  
scattering objects with elliptical geometry.  The fast computation and the information from the results of the sampling  method help
speed up the training of the network.  Our numerical study shows that  the network  is able to
predict  some scattering objects with geometries different from those of the training and testing data sets. The combination leads 
to a significant improvement in the reconstruction results initially obtained by the sampling method.
 We also demonstrate that the combined method is also able to invert some experimental data
from the Fresnel Institute (\cite{Belke2001}) without any additional transfer learning.

The paper is organized as follows. The orthogonality sampling type method and its analysis are presented in Section 2. Section 3 is dedicated to a description of the deep neural network 
used to combine with the sampling method. The numerical study for simulated data and real data is presented in Sections 4 and 5, respectively.

\section{An orthogonality sampling type method}
In this section we develop the an orthogonality sampling type method for the inverse scattering problem. For $x, y \in \R^3$ and $x \neq y$, we
denote by $\Phi(x,y)$ the free-space Green's function of the scattering problem~\eqref{Helm}--\eqref{radiation} which is given by 
\begin{equation} 
\label{green}
 \Phi(x,y)= 
\begin{cases}
\frac{i}{4}H^{(1)}_0(k|x-y|), & n=2 , \\ 
\frac{\exp(ik|x-y|)}{4\pi|x-y|}, & n=3.
\end{cases}
\end{equation}
 It is  well known that  problem~\eqref{Helm}--\eqref{radiation} is equivalent to the Lippmann-Schwinger equation (see, e.g.~\cite{Colto2013})
$$
u_\sc(x) = k^2\int_D \Phi(x,y) \eta(y) u(y)\dd y, \quad x \in \R^3.
$$
The following Helmholtz integral representation is important to the  analysis of the sampling method. 
\begin{lemma}
Assume that $W$ is  bounded Lipschitz  domain in $\R^n$    such that $\R^n \setminus \overline{W}$ is connected. Let $\nu$ be a unit outward normal vector on $\partial W$. 
For any function $ w \in H^2(W)$ we
have
\begin{align}
\label{id}
w(x) = \int_{\partial W} \left(\frac{\partial w(y)}{\partial \nu} \Phi(x,y) - w(y) \frac{\Phi(x,y)}{\partial\nu(y)} \right) \dd s(y) -  \int_{ W} (\Delta w + k^2 w) \Phi(x,y) \dd y.
\end{align}
\end{lemma}
\begin{proof}
We refer to~\cite[Chapter  2]{Colto2013} for a proof of the lemma.
\end{proof} 
We also need to the Funk-Hecke formula for the analysis in this section.
\begin{lemma}
Let $Y_m$ be the spherical harmonics of order $m$ and $j_m$ be the spherical Bessel functions of first kind and order $m$. We have
\begin{align}
\label{funk}
\int_{\S^{n-1}} e^{-ikx \cdot \hat{z}} Y_n(\hat{z}) ds(\hat{z}) = \frac{4 \pi}{i^m}j_m(k|x|)Y_m\left(\frac{x}{|x|}\right).
\end{align}
\end{lemma}
\begin{proof}
A proof of this lemma can be found in~\cite[Chapter 2]{Colto2013}. 
\end{proof} 
Now we define the imaging functional as
\begin{align}
\label{I2}
\mathcal{I}(z) :=  \left | \int_{\partial \Omega} \left(\frac{\partial \Im \Phi(x,z)}{\partial \nu(x)} u_\sc(x) - \Im \Phi(x,z) \frac{\partial u_\sc(x)}{\partial \nu(x)} \right) \dd s(x)\right |^\rho,
\end{align}
where $\rho = 1$ or  $\rho = 2$. 

This functional $\mathcal{I}(z)$ aims to determine $D$ in  $\Omega$ and in the numerical implementation 
we evaluate $\mathcal{I}(z)$ for a finite set of sampling points  $z$. We expect that $\mathcal{I}(z)$ is relatively large
as $z \in D$ and that it is small as $z$ is outside $D$.

The   behavior of $\mathcal{I}(z)$ is analyzed in the following theorem.
\begin{theorem}
The imaging functional $\mathcal{I}$ satisfies 
$$
\mathcal{I}(z) = \left |k^2 \int_D  \Im \Phi(y,z) \eta(y) u(y) \dd y\right |^\rho.
$$
\end{theorem}
\begin{proof}
From the Lippmann-Schwinger equation we obtain that
$$
\frac{\partial u_\sc(x)}{\partial \nu(x)} = k^2\int_D \frac{\partial \Phi(x,y)}{\partial \nu(x)} \eta(y)u(y) \dd y.
$$ 
Therefore, a substitution  leads to 
\begin{align}
& \int_{\partial \Omega} \left(\frac{\partial \Im \Phi(x,z)}{\partial \nu(x)} u_\sc(x) - \Im \Phi(x,z) \frac{\partial u_\sc(x)}{\partial \nu(x)} \right) \dd s(x) \nonumber \\
 & =  \int_{\partial\Omega} \left (\frac{\partial \Im \Phi(x,z)}{\partial \nu(x)} k^2\int_D  \Phi(x,y) \eta(y)u(y) \dd y - \Im \Phi(x,z) k^2  \int_\Omega \frac{\partial  \Phi(x,y)}{\partial \nu(x)}\eta(y)u(y) \dd y \right) \dd s(x) \nonumber \\
\label{int1}
 &= k^2 \int_D \int_{\partial \Omega}  \left(\frac{\partial \Im \Phi(x,z)}{\partial \nu(x)} \Phi(y,x) - \Im \Phi(x,z) \frac{\partial \Phi(y,x)}{\partial \nu(x)} \right)\dd s(x) \, \eta(y) u(y) \dd y.
\end{align}
Since $\Delta \Im \Phi(z,y) + k^2 \Im \Phi(z,y) = 0$ for  all $z, y \in \mathbb{R}^n$ and  $\Im \Phi(z,y)$ is regular function,  the Helmholtz integral representation~\eqref{id} implies that
$$
  \int_{\partial \Omega}  \left(\frac{\partial \Im \Phi(x,z)}{\partial \nu(x)} \Phi(y,x) - \Im \Phi(x,z) \frac{\partial \Phi(y,x)}{\partial \nu(x)} \right)\dd s(x) = \Im \Phi(y,z).
$$
Therefore, substituting this identity in~\eqref{int1} implies
$$
 \int_{\partial \Omega} \left(\frac{\partial \Im \Phi(x,z)}{\partial \nu(x)} u_\sc(x) - \Im \Phi(x,z) \frac{\partial u_\sc(x)}{\partial \nu(x)} \right) \dd s(x)  = k^2\int_D\Im \Phi(y,z) \eta(y) u(y) \dd y.
 $$
Now the proof can be completed  by substituting this equation in~\eqref{I2}.
\end{proof}
We know that 
\begin{equation} 
\Im \Phi(y,z)= 
\begin{cases}
J_0(k|z-y|), & n=2 , \\ 
\frac{k}{4\pi}j_0(k|z-y|), & n=3.
\end{cases}
\label{Im}
\end{equation}
Since $J_0(k|y-z|)$ and $j_0(k|y-z|)$ peak when $k|y-z| = 0$, we expect from Theorem 1 that $\mathcal{I}(z)$ takes larger values for $z \in D$ 
and  much smaller values outside $D$. Further, from the asymptotic behavior of $J_0(k|y-z|)$ and $j_0(k|y-z|)$ as $|y-z| \to \infty$ we can easily estimate, for $z \notin D$, that
$$
\mathcal{I}(z) = O\left(\frac{1}{\mathrm{dist}(z,D)^{\rho(n-1)/2}} \right) \text{ as } \mathrm{dist}(z,D) \to \infty.
$$
This incomplete analysis is common for imaging functionals of orthogonality sampling methods (see, e.g., ~\cite{Potth2010, Ito2012, Ito2013, Akinc2016, Nguye2019, Ahn2020}). 
To  our knowledge a complete analysis for orthogonality sampling methods is still an open problem.

We now  show a relation between $\mathcal{I}(z)$  and the imaging functional suggested by Potthast in~\cite{Potth2010}. To this end, 
we recall that  the scattered field $u_\sc$
has the asymptotic behavior
$$
u_\sc(x) = \frac{e^{ik|x|}}{|x|^{(n-1)/2}}(u^\infty(\hat{x}) + O(1/|x|)), \quad |x| \to \infty,
$$
for all $\hat{x} \in \S^{n-1}$. The function $u^\infty(\hat{x})$ is called the scattering amplitude or the far-field pattern of $u_\sc(x)$. The original imaging functional of the orthogonality sampling method suggested 
in~\cite{Potth2010} is given by
$$
\mathcal{I}_{\mathrm{OSM}}(z) = \left|\int_{\S^{n-1}} e^{ikz\cdot \hat{x}} u^\infty(\hat{x}) \dd s(\hat{x})\right|.
$$
\begin{theorem}
\label{relation}
The two imaging functionals are related by
$$
\mathcal{I}(z) = \gamma  \left| \mathcal{I}_{\mathrm{OSM}}(z)  \right|^\rho 
$$
where
 \begin{equation*} 
\gamma = 
\begin{cases}
\left(\frac{\sqrt{\pi}e^{i\pi/4}}{\sqrt{2 k}} \right)^\rho, & n= 2 , \\ 
\left(\frac{4\pi}{k} \right)^\rho, & n= 3.
\end{cases}
\end{equation*}
\end{theorem}
\begin{proof}
From the Lippmann-Schwinger  equation for $u_\sc$ it is  known that  its scattering amplitude is given by 
$$
u^\infty(\hat{x}) = \alpha k^2\int_D e^{-iky\cdot \hat{x}} \eta(y) u(y)\dd y
$$
with the constant
 \begin{equation*} 
\alpha = 
\begin{cases}
\frac{e^{i\pi/4}}{\sqrt{8\pi k}}, & n= 2 , \\ 
\frac{1}{4\pi}, & n= 3.
\end{cases}
\end{equation*}
Substituting $u^\infty(\hat{x})$ in $\mathcal{I}_{\mathrm{OSM}}(z)$ implies 
\begin{equation}
\label{IOSM}
\mathcal{I}_{\mathrm{OSM}}(z) =  \left|\alpha k^2\int_D \int_{\S^{n-1}}  e^{ik\hat{x}\cdot (z-y)} \dd s(\hat{x}) \eta(y)u(y) \dd y \right|.
\end{equation}
Now using~\eqref{Im} and  the Funk-Hecke formula~\eqref{funk} we obtain that
\begin{equation} 
\int_{\S^{n-1}} e^{ik\hat{x}\cdot (z-y)} \dd s(\hat{x})= 
\begin{cases}
2\pi \Im \Phi(z,y), &n=2 , \\ 
\frac{(4\pi)^2}{k}  \Im \Phi(z,y), & n=3.
\end{cases}
\label{FH}
\end{equation}
Substituting~\eqref{FH} in~\eqref{IOSM} and comparing with Theorem 1 we complete the proof.
\end{proof}

Theorem 2 also implies that $\mathcal{I}(z)$ is stable against noise in the data following the stability of $\mathcal{I}_{\mathrm{OSM}}(z)$~\cite{Potth2010}.

\begin{remark}
If $\partial \Omega$ is the boundary of some ball with large radius, we can approximate  $\partial u_\sc/\partial \nu \approx ik u_\sc$ using the  radiation condition. 
Then   $\mathcal{I}(z)$ can be modified  to handle  the far-field  data $u_\sc(x)$ as follows
\begin{align}
\label{I2far}
\mathcal{I}(z) =  \left | \int_{\partial \Omega} \left(\frac{\partial \Im \Phi(x,z)}{\partial \nu(x)} u_\sc(x) - ik \Im \Phi(x,z)  u_\sc(x) \right) \dd s(x)\right |^\rho.
\end{align}
\end{remark}


\section{The deep neural network}
In this section, we use the sampling method along with a deep neural network of U-Net-Xception style to form a combined method to solve the inverse problem. Numerical simulation was done for both simulated and experimental scattering data.

\subsection{U-Net Xception style}
Xception is an efficient architecture which relies on two main points: Depthwise Separable Convolution and Shortcuts between Convolution blocks. Xception architecture is an alternative to classical convolution layers, which performs more efficiently in terms of computation time and accuracy. See details of Xception architecture in~\cite{Choll2016}.

\begin{figure}[tb]
\centering
\includegraphics[width = 0.85\textwidth]{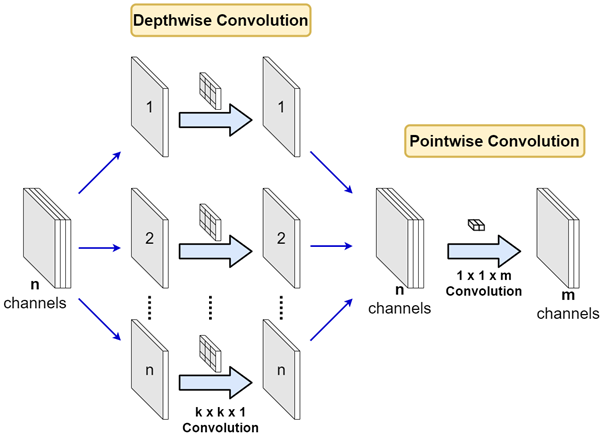}
\caption{Architecture of Depthwise Separable Convolution.}
\label{xcep}
\end{figure}

Firstly, Depthwise Separable Convolution includes two subnetworks: Depthwise Convolution and Pointwise Convolution. Compared to conventional convolution layers, Depthwise Separable Convolution does not need to compute convolution operations across all channels, which will eliminate the number of parameters in the model and makes it simpler. The overall architecture of Depthwise Separable Convolution is described in Figure \ref{xcep}.

\begin{figure}[tb]
\centering
\includegraphics[width = 0.5\textwidth]{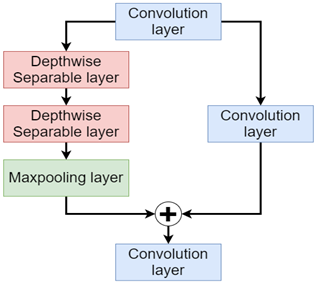}
\caption{Residual connections within Xception networks.}
\label{connection}
\end{figure}

Secondly, there are residual (shortcut/skip) connections in Xception networks. Many state-of-the-art Deep neural networks (Resnet, YOLOv3, MobileNetV2, Cascading Residual Network) have shown that accuracy is much higher when residual connections are used, see \cite{He2015, Redmo2018, Sandl2019, Ahn2018}. The pattern of residual connection in U-Net Xception architecture is shown in Figure \ref{connection}. 

\subsection{Training the network}

Given some scattering object and an incident wave we use the spectral method studied in~\cite{Lechl2014} to solve the direct  problem and generate the associated scattering data. 
The image of the  scattering object is called a true image. Applying the sampling method  to the scattering data we 
obtain a computed image of the scattering object. We call this image a preliminary image.  

\begin{figure}[tb]
\centering
\begin{subfigure}{0.3\textwidth}
\centering
\includegraphics[width=0.85\textwidth]{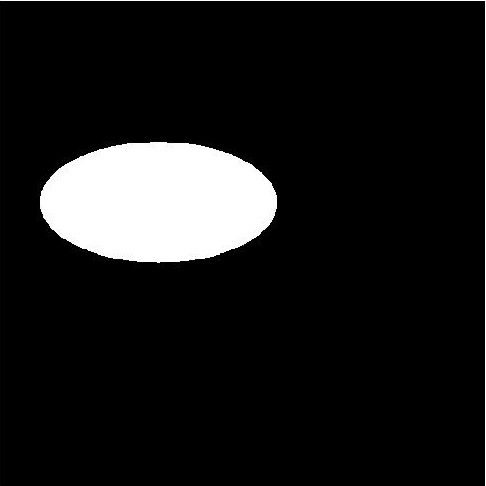}
\end{subfigure}
\begin{subfigure}{0.3\textwidth}
\centering
\includegraphics[width=0.85\textwidth]{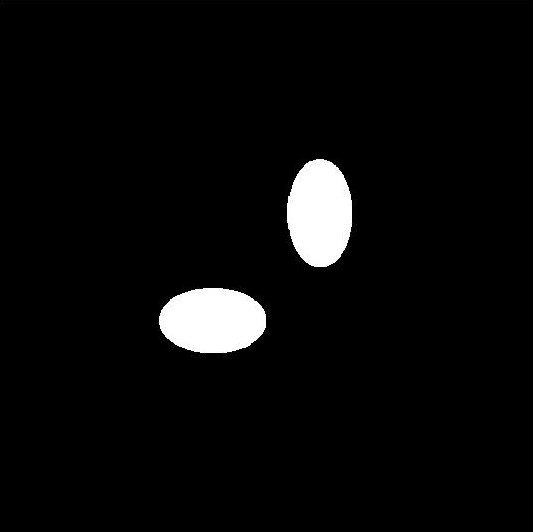}
\end{subfigure}
\begin{subfigure}{0.3\textwidth}
\centering
\includegraphics[width=0.85\textwidth]{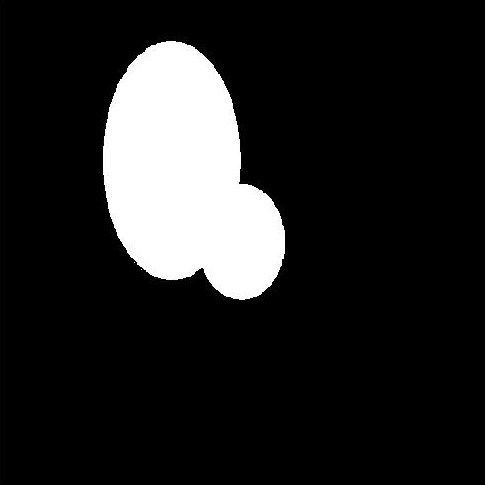}
\end{subfigure}

\vspace{2mm}

\begin{subfigure}{0.3\textwidth}
\centering
\includegraphics[width=0.85\textwidth]{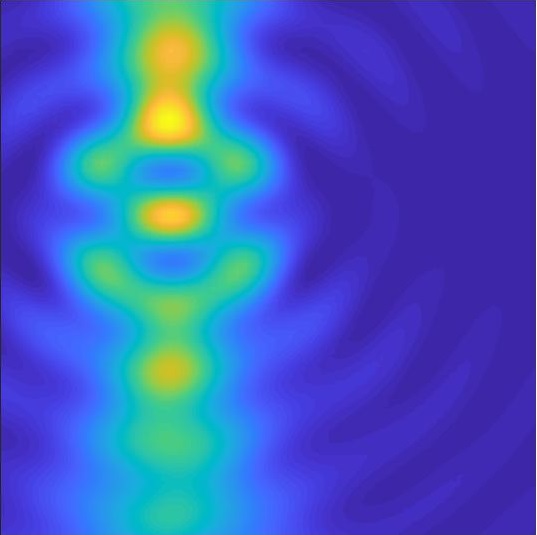}
\end{subfigure}
\begin{subfigure}{0.3\textwidth}
\centering
\includegraphics[width=0.85\textwidth]{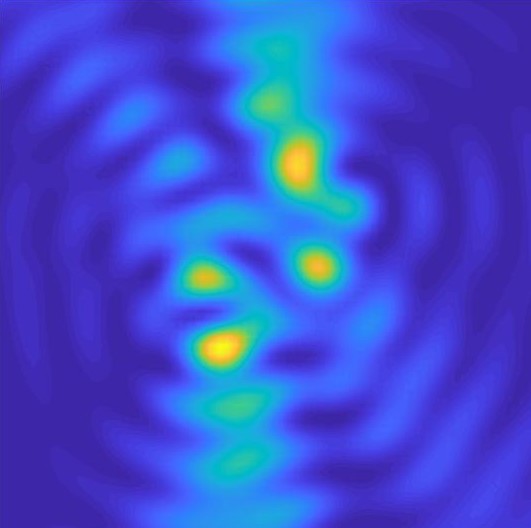}
\end{subfigure}
\begin{subfigure}{0.3\textwidth}
\centering
\includegraphics[width=0.85\textwidth]{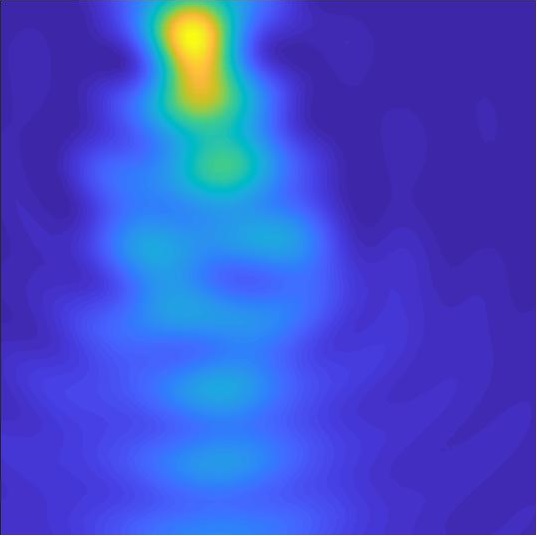}
\end{subfigure}
\caption{Examples of three training data pairs: true images (top), preliminary images (bottom).}
\label{example}
\end{figure}

The training data set includes $60,000$ pairs of images, each pair consists of a true image and a preliminary image. The true images are those of either one or two ellipses whose center and radii are randomly generated within the square domain $(-2,2)^2$. More precisely, if the object consists of only one ellipse, the center is randomly generated within the square $[-0.8,0.8]^2$ and the radii are randomly generated within the interval $[0.1,1]$; if the object consists of two ellipses, the first one is generated in the same manner as before, while the second one has its center randomly generated within the square $[-1,1]^2$ and radii randomly generated within the interval $[0.1,0.5]$. Note that the two ellipses are allowed to overlap each other. The wave number associated with the data set is $k = 6$, which means the wavelength 
is about 1.  As for the incident wave, we use $u^{in}(x) = e^{ik(x_1\cos\theta+x_2\sin\theta)}$ where $\theta = 90^o$ for 30,000 images and $\theta = 45^o$ for the other 30,000 images. 
The Cauchy scattering data are computed on the circle of radius 100. We did not add any artificial noise to the scattering data for the training and testing of the network.  See Figure \ref{example} for some examples of training data pairs.

\begin{figure}[tb]
\centering
\includegraphics[width = \textwidth]{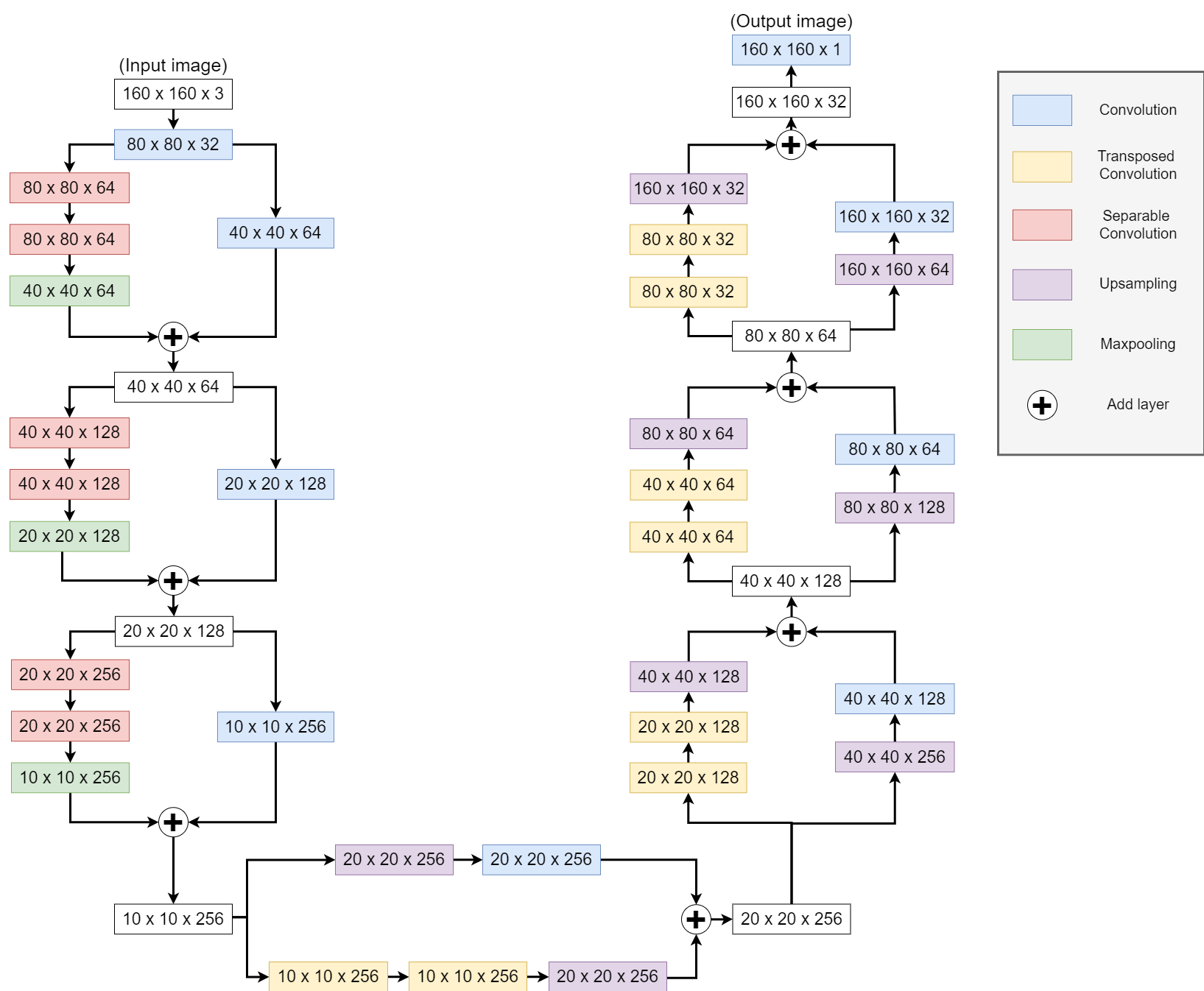}
\caption{Architecture of the neural network.}
\label{architecture}
\end{figure}

All images are uniformly partitioned into $160 \times 160$ pixels. Let $I = \{1,2,\dots,160^2\}$. The neural network is initialized with random parameters. Using a preliminary image as input, the network produces the corresponding output $\hat{y} = \{\hat{y}_i\}_{i\in I}$ that will then be used along with the true image $y = \{y_i\}_{i\in I}$ to calculate the loss function
$$
L = -\sum_{i\in I}(y_i\log \hat{y}_i + (1-y_i)\log(1-\hat{y}_i)).
$$
After that, the parameters within the convolution layers are tuned using Adam method (see \cite{Kingm2014}) in order to minimize the loss function. This minimization process is done for one batch of size $32$ at a time, with the loss function for each batch being the sum of those of individual training pairs in the batch. The initial learning rate is $10^{-2}$, and will be divided by $10$ if the validation loss does not decrease or validation accuracy does not increase after $5$ epochs until it reaches $10^{-6}$. The neural network's architecture can be briefly described by the diagram in Figure \ref{architecture}.

\begin{figure}[tb]
\centering
\begin{subfigure}{0.49\textwidth}
\centering
\includegraphics[width=\textwidth]{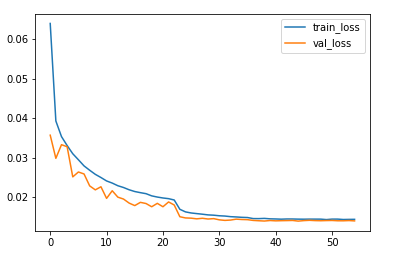}
\caption{Loss}
\end{subfigure}
\begin{subfigure}{0.49\textwidth}
\centering
\includegraphics[width=\textwidth]{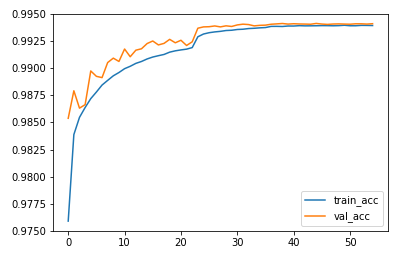}
\caption{Accuracy}
\end{subfigure}
\caption{Training and validation results through 50 epochs.}
\label{test_res}
\end{figure}

Out of 60,000 pairs of images, 50,000 were used to train the network, 5,000 were used for validation and 5,000 are used for testing. In Figure \ref{test_res}, we show the result of training loss versus validation loss and training accuracy versus validation accuracy. The loss value at each epoch is the average value of $L$ over all paired images, whereas the accuracy value is the ratio of the number of pixels in the predicted image that match their counterparts in the true image to the total number of pixels. Those values for validation were computed over $5,000$ images in the validation set, and those values for training were computed over $50,000$ images in the training set. 

\begin{figure}[tb]
\centering
\begin{subfigure}{0.3\textwidth}
\centering
\includegraphics[width=0.85\textwidth]{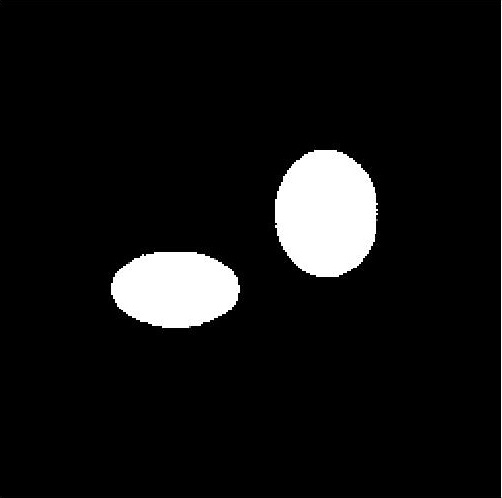}
\caption{Original image}
\end{subfigure}
\begin{subfigure}{0.3\textwidth}
\centering
\includegraphics[width=0.85\textwidth]{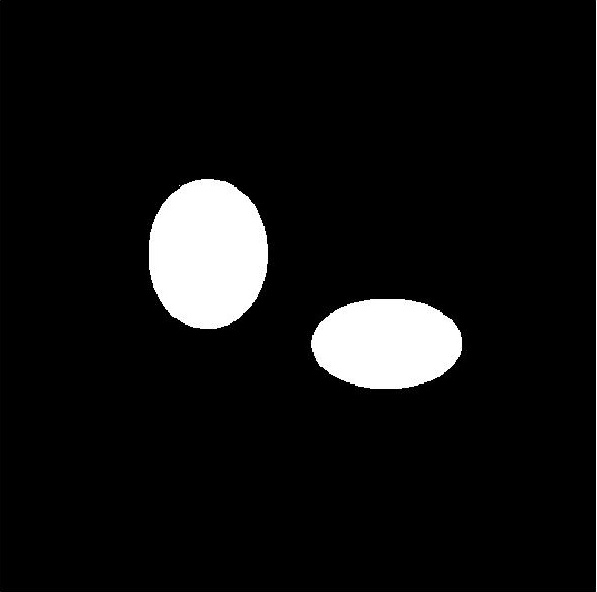}
\caption{Horizontal flip}
\end{subfigure}
\begin{subfigure}{0.3\textwidth}
\centering
\includegraphics[width=0.85\textwidth]{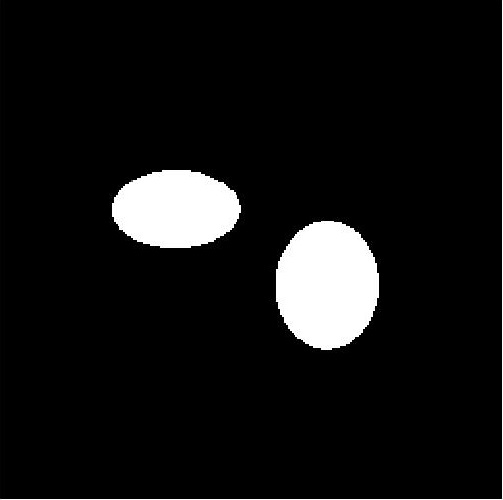}
\caption{Vertical flip}
\end{subfigure}

\vspace{2mm}

\begin{subfigure}{0.3\textwidth}
\centering
\includegraphics[width=0.85\textwidth]{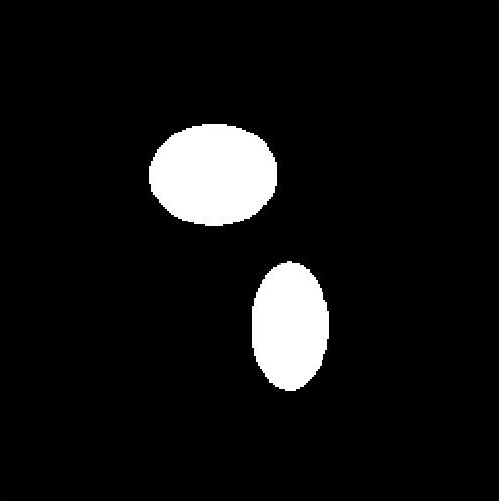}
\caption{$90^o$ rotation}
\end{subfigure}
\begin{subfigure}{0.3\textwidth}
\centering
\includegraphics[width=0.85\textwidth]{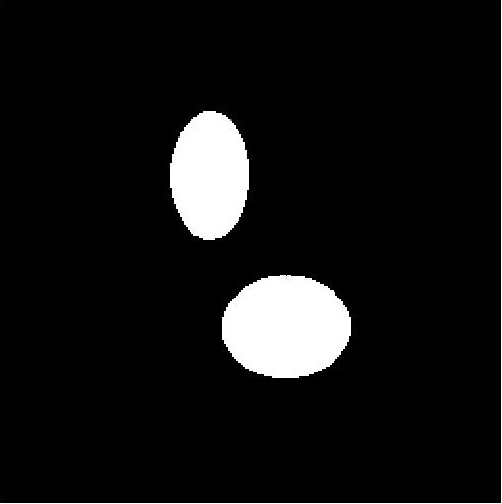}
\caption{$-90^o$ rotation}
\end{subfigure}
\begin{subfigure}{0.3\textwidth}
\centering
\includegraphics[width=0.85\textwidth]{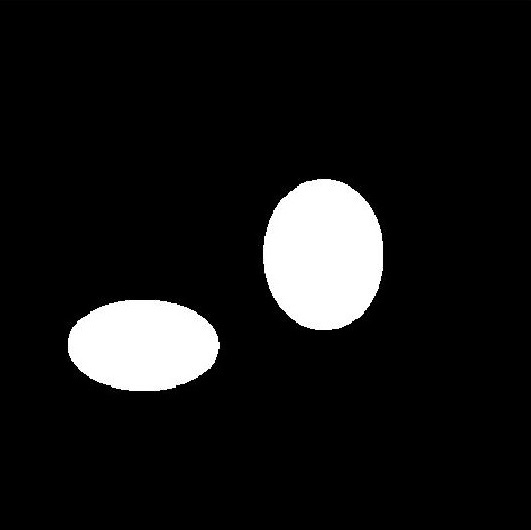}
\caption{Zoom}
\end{subfigure}
\caption{An example of data augmentation}
\label{aug}
\end{figure}

In order to enrich the training data set, we used \textit{data augmentation}. More specifically, we produced five copies of each training pair by applying five types of transformations: horizontal flip, vertical flip, $90^o$ rotation, $-90^o$ rotation and zoom. Zoom transformations were done by first upscaling the image by $5\%$, then using a random frame of the same size as the original image to crop it. See Figure \ref{aug} for an example of such augmentations.

Our choice for some hyperparameters was based on direct comparison to other choices in terms of validation loss and accuracy. For example, we trained the model with $10,000$ pairs of images through $5$ epochs with different input image sizes, then we chose $160\times 160\times 3$ because it gave the best values for loss and accuracy. We did a similar experiment for two different optimization methods, which are Adam method and stochastic gradient descent with momentum (SGD), and Adam method performed better for this problem. Detailed results for the comparisons are given in Table \ref{iisize} and \ref{omethod}. As for the other hyperparameters such as the number of hidden layers and units, we followed the same choice as in~\cite{chollet2015keras}. With this choice, we were able to get very good loss and accuracy values, see Figure~\ref{test_res}. Therefore, no adjustment was needed.

\begin{table}[tb]
\centering
\renewcommand{\arraystretch}{2}
\begin{tabular}{|c|c|c|}
\hline
\phantom{a} Input image size \phantom{a} & \phantom{a} Validation accuracy \phantom{a} &  \phantom{aa} Validation loss  \phantom{aa} \\
\hline
$32\times 32\times 3$ & $0.9718$ & $0.0677$ \\
\hline
$64\times64\times3$ & $0.9762$ & $0.0579$ \\
\hline
$128\times128\times3$ & $0.9728$ & $0.0713$ \\
\hline
$\mathbf{160\times160\times3}$ & $\mathbf{0.9822}$ & $\mathbf{0.0452}$  \\
\hline
$192\times192\times3$ & $0.9749$ & $0.0708$ \\
\hline
\end{tabular}
\caption{Input image size comparison}
\label{iisize}
\end{table}

\begin{table}[tb]
\centering
\renewcommand{\arraystretch}{2}
\begin{tabular}{|l|c|c|}
\hline
\phantom{a} Optimization method \phantom{a}  & \phantom{a} Validation accuracy \phantom{a} &  \phantom{aa} Validation loss  \phantom{aa} \\
\hline
$\mathrm{SGD}$ & $0.9751$ & $0.0622$ \\
\hline
\textbf{Adam} & $\mathbf{0.9822}$ & $\mathbf{0.0452}$ \\
\hline
\end{tabular}
\caption{Optimization method comparison}
\label{omethod}
\end{table}

We also compared the performance of the U-Net Xception style with a modified U-Net  in \cite{daous2019tensorflow}. The result in Table \ref{model} shows that U-Net Xception style works better for our problem.

\begin{table}[tb]
\centering
\renewcommand{\arraystretch}{2}
\begin{tabular}{|c|c|c|}
\hline
Model  & \phantom{a} Validation accuracy \phantom{a} &  \phantom{aa} Validation loss  \phantom{aa}  \\
\hline
Modified U-Net \phantom{aaaaaaaaa} & $0.9876$ & $0.0299$ \\
\hline
\textbf{U-Net Xception style} & $\mathbf{0.9932}$ & $\mathbf{0.0175}$ \\
\hline
\end{tabular}
\caption{Model comparison}
\label{model}
\end{table}

\section{Numerical study for simulated  data}
In this section we use the combined method to reconstruct objects from simulated scattering data. We consider the types of objects varying from those similar to the training data set to those entirely different. We also tested different wave numbers and incident wave directions.

\subsection{Implementation of the sampling method}
In our numerical study the sampling domain is chosen as $[-2,2]^2$.  The sampling domain contains the unknown 
scatterer  $D$ and is uniformly discretized by 64 sampling points in each direction. 
The imaging functional $\mathcal{I}(z)$ in~\eqref{I2} is evaluated at each sampling point $z$ of the 
sampling domain. As mentioned in the section of training the neural network  
the  data $(u_\sc, \partial u_\sc/ \partial \nu)$ are given at 32 uniformly distributed points  
on the circle of  radius 100.  With the data given, the evaluation of $\mathcal{I}(z)$  is simple since  
after a change of variables using polar coordinates we  numerically evaluate a 
single integral  using the rectangular rule. In the numerical simulation  the imaging functional is 
normalized by dividing by its (numerical) maximal value. 
 Note that we make sure that the sampling point $z$ and the data point $x$ are away from 
 each other so that  the  integrands are all regular functions. Recall that the scattering data 
 are measured on $\partial \Omega$ meaning that we 
 should look for the unknown scatterer $D$ in a sampling domain that is  strictly contained in $\Omega$.

\subsection{Reconstruction of objects in the testing data set}
In this part, we present some reconstruction results of the combined method  for objects in the testing data set.  These objects were generated in a way that is similar to that of the training data sets, meaning the objects consist of one or two ellipses. The wave number $k$ and incident direction $\theta$ were those the neural network was trained with. Moreover, the scattering data generated by solving the direct problem was not perturbed by noise.

\begin{figure}[tb]
\centering
\begin{subfigure}{0.3\textwidth}
\centering
\includegraphics[width=0.85\textwidth]{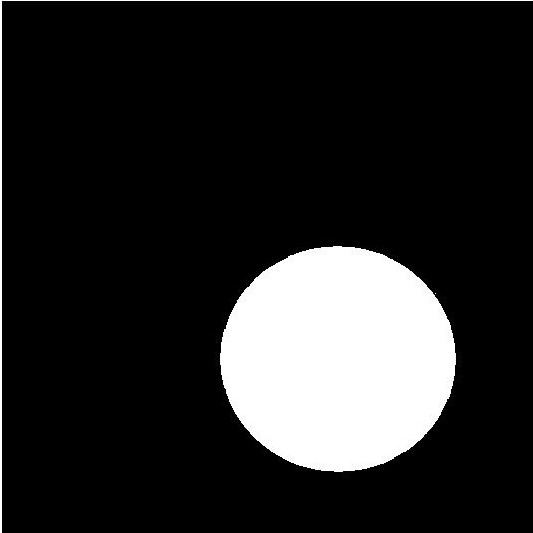}
\end{subfigure}
\begin{subfigure}{0.3\textwidth}
\centering
\includegraphics[width=0.85\textwidth]{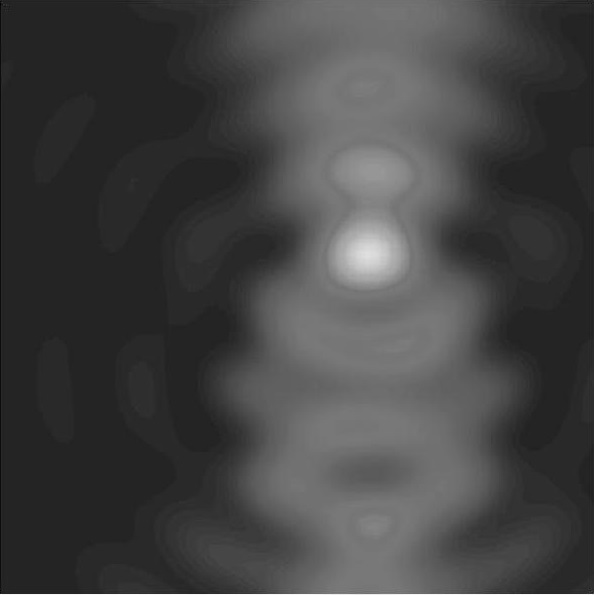}
\end{subfigure}
\begin{subfigure}{0.3\textwidth}
\centering
\includegraphics[width=0.85\textwidth]{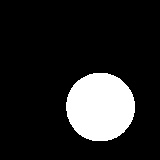}
\end{subfigure}

\vspace{3mm}

\begin{subfigure}{0.3\textwidth}
\centering
\includegraphics[width=0.85\textwidth]{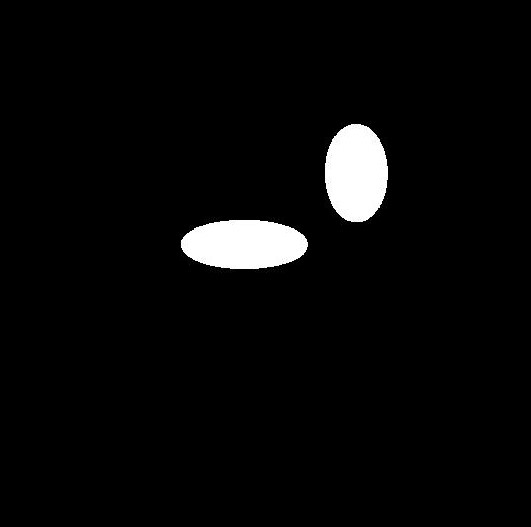}
\caption{True objects}
\end{subfigure}
\begin{subfigure}{0.3\textwidth}
\centering
\includegraphics[width=0.85\textwidth]{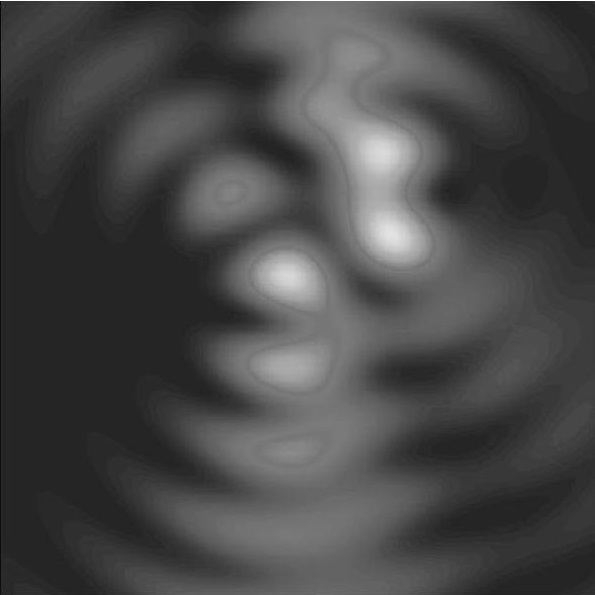}
\caption{Sampling method}
\end{subfigure}
\begin{subfigure}{0.3\textwidth}
\centering
\includegraphics[width=0.85\textwidth]{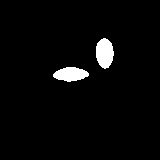}
\caption{Combined method}
\end{subfigure}
\caption{Reconstruction of objects in the testing data set with $k = 6$ and $\theta = 90^o$: one disk (top) and two ellipses of similar sizes (bottom).}
\label{fig7}
\end{figure}

\begin{figure}[tb]
\centering
\begin{subfigure}{0.3\textwidth}
\centering
\includegraphics[width=0.85\textwidth]{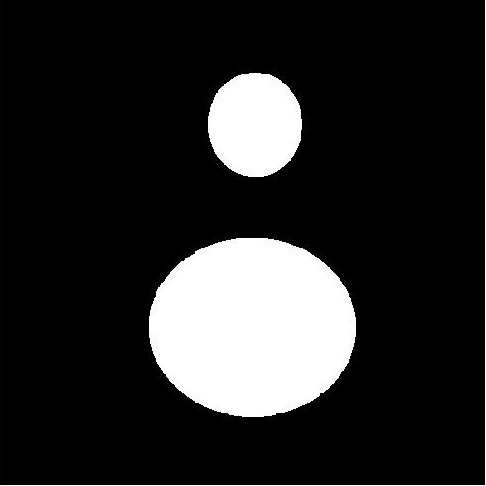}
\end{subfigure}
\begin{subfigure}{0.3\textwidth}
\centering
\includegraphics[width=0.85\textwidth]{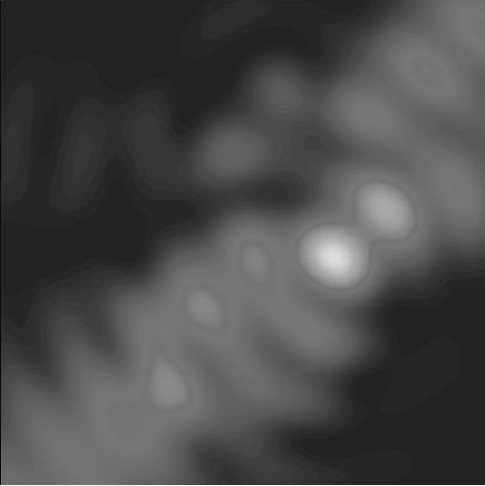}
\end{subfigure}
\begin{subfigure}{0.3\textwidth}
\centering
\includegraphics[width=0.85\textwidth]{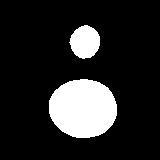}
\end{subfigure}

\vspace{3mm}

\begin{subfigure}{0.3\textwidth}
\centering
\includegraphics[width=0.85\textwidth]{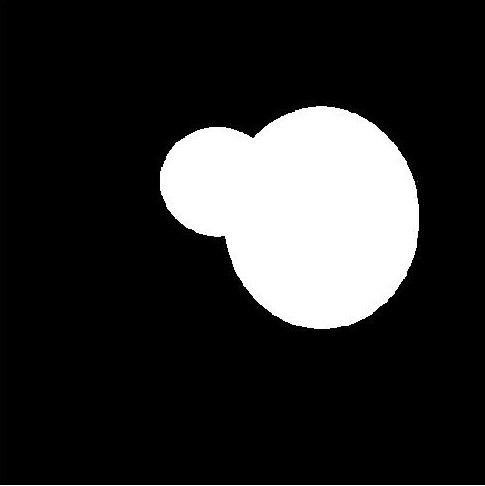}
\caption{True objects}
\end{subfigure}
\begin{subfigure}{0.3\textwidth}
\centering
\includegraphics[width=0.85\textwidth]{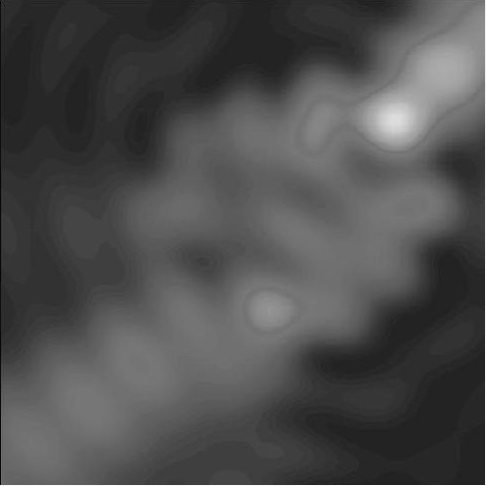}
\caption{Sampling method}
\end{subfigure}
\begin{subfigure}{0.3\textwidth}
\centering
\includegraphics[width=0.85\textwidth]{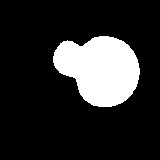}
\caption{Combined method}
\end{subfigure}
\caption{Reconstruction of objects in the  testing data set with $k = 6$ and $\theta = 45^o$:  two ellipses of different sizes (top) and two overlapping ellipses (bottom).}
\label{fig8}
\end{figure}

 In Figures~\ref{fig7}--\ref{fig8}, four test objects are presented: one large disk, two ellipses of similar sizes, two ellipses of different sizes and two overlapping ellipses. We can see that the method was able to reconstruct  very well the geometry of these test objects. As it was done for the validation and training data sets in Figure~\ref{test_res}, the  accuracy  for  5000 pairs of images in the testing data set  is also computed as $99.4\%$. This once again confirms that the training process was effective as indicated by the training and validation accuracy in Figure \ref{test_res}.

%

\subsection{Reconstruction of objects with no elliptical shapes}
\begin{figure}[tb]
\centering
\begin{subfigure}{0.3\textwidth}
\centering
\includegraphics[width=0.85\textwidth]{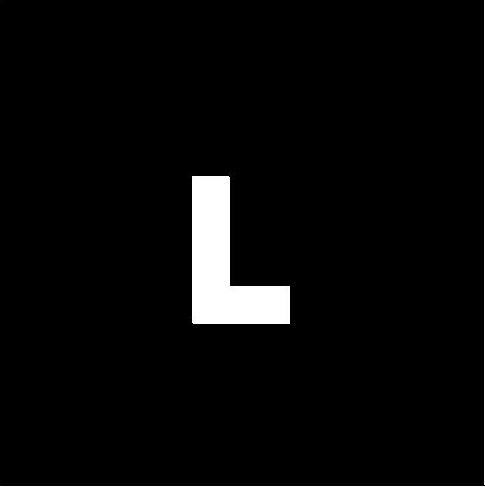}
\caption{True object}
\end{subfigure}
\begin{subfigure}{0.3\textwidth}
\centering
\includegraphics[width=0.85\textwidth]{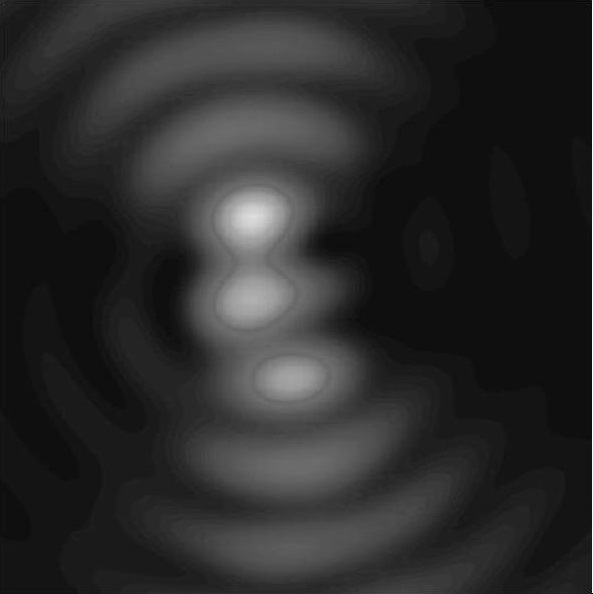}
\caption{Sampling, $5\%$ noise}
\end{subfigure}
\begin{subfigure}{0.3\textwidth}
\centering
\includegraphics[width=0.85\textwidth]{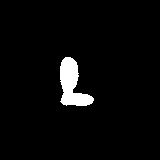}
\caption{Combined, $5\%$ noise}
\end{subfigure}

\vspace{2mm}

\begin{subfigure}{0.3\textwidth}
\centering
\includegraphics[width=0.85\textwidth]{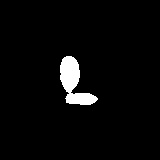}
\caption{Combined, $7\%$ noise}
\end{subfigure}
\begin{subfigure}{0.3\textwidth}
\centering
\includegraphics[width=0.85\textwidth]{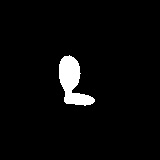}
\caption{Combined, $10\%$ noise}
\end{subfigure}
\begin{subfigure}{0.3\textwidth}
\centering
\includegraphics[width=0.85\textwidth]{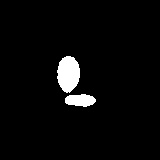}
\caption{Combined, $15\%$ noise}
\end{subfigure}

\caption{Reconstruction of an L-shaped object with $k = 6$, $\theta = 90^o$, and different noise levels added to the scattering data: $5\%$, $7\%$, $10\%$ and $15\%$.}
\label{L}
\end{figure}

In this part, we test the method for other types of objects with noisy scattering data.  We added $5\%, 7\%, 10\%$ and $15\%$ of artificial noise to the scattering data obtained from solving the direct problem with wave number $k=6$. Recall that the inverse problem with one incident wave is extremely ill-posed, and there were no artificial noise added to the scattering data in the training process. The tests with noisy data for non-elliptical objects are to show that the network can generalize its training to more realistic situations that are unseen from the training process. We first used the same incident direction as in training, $\theta = 90^o$, to reconstruct  an L-shaped object, see Figure \ref{L}.
The combined method was able to provide reasonable reconstructions for the  L-shaped object with  different levels of noise.

\begin{figure}[tb]
\centering
\begin{subfigure}{0.3\textwidth}
\centering
\includegraphics[width=0.85\textwidth]{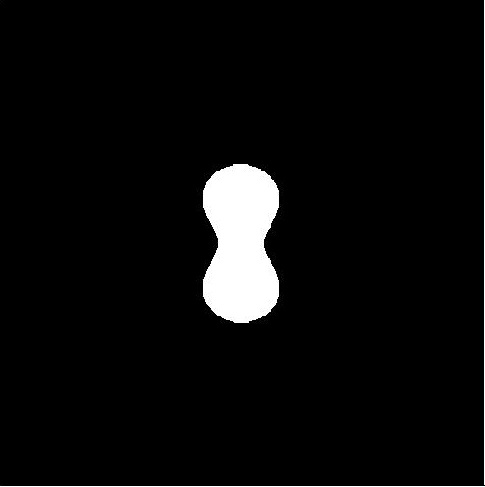}
\caption{True object}
\end{subfigure}
\begin{subfigure}{0.3\textwidth}
\centering
\includegraphics[width=0.85\textwidth]{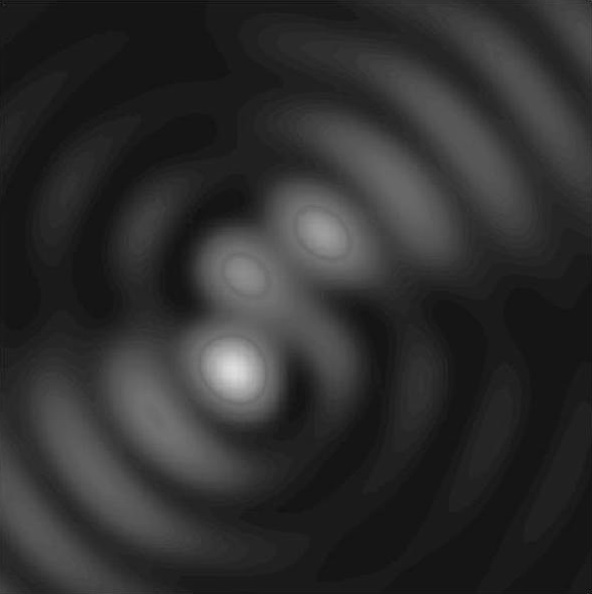}
\caption{Sampling, $5\%$ noise}
\end{subfigure}
\begin{subfigure}{0.3\textwidth}
\centering
\includegraphics[width=0.85\textwidth]{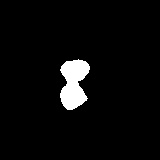}
\caption{Combined, $5\%$ noise}
\end{subfigure}

\vspace{2mm}

\begin{subfigure}{0.3\textwidth}
\centering
\includegraphics[width=0.85\textwidth]{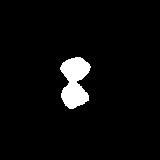}
\caption{Combined, $7\%$ noise}
\end{subfigure}
\begin{subfigure}{0.3\textwidth}
\centering
\includegraphics[width=0.85\textwidth]{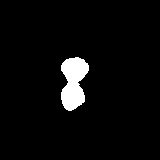}
\caption{Combined, $10\%$ noise}
\end{subfigure}
\begin{subfigure}{0.3\textwidth}
\centering
\includegraphics[width=0.85\textwidth]{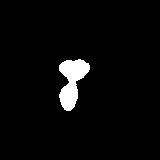}
\caption{Combined, $15\%$ noise}
\end{subfigure}
\caption{Reconstruction of a peanut-shaped object with $k = 6$, $\theta = 220^o$, and different noise levels added to the scattering data: $5\%$, $7\%$, $10\%$ and $15\%$.}
\label{pea}
\end{figure}

\begin{figure}[tb]
\centering
\begin{subfigure}{0.3\textwidth}
\centering
\includegraphics[width=0.85\textwidth]{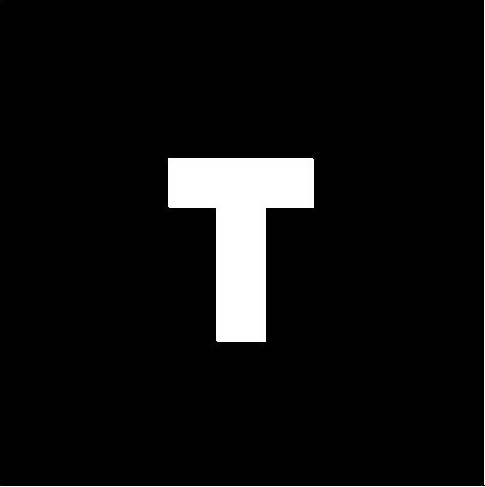}
\caption{True object  }
\end{subfigure}
\begin{subfigure}{0.3\textwidth}
\centering
\includegraphics[width=0.85\textwidth]{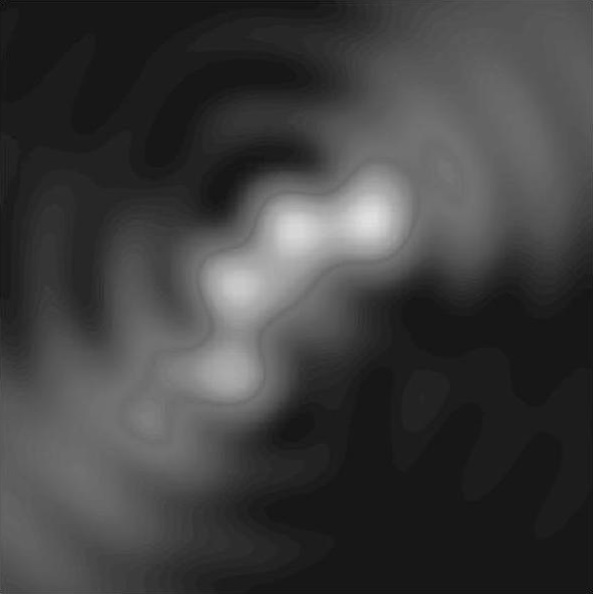}
\caption{Sampling, $5\%$ noise}
\end{subfigure}
\begin{subfigure}{0.3\textwidth}
\centering
\includegraphics[width=0.85\textwidth]{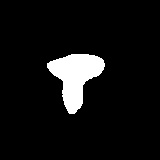}
\caption{Combined, $5\%$ noise}
\end{subfigure}

\vspace{2mm}

\begin{subfigure}{0.3\textwidth}
\centering
\includegraphics[width=0.85\textwidth]{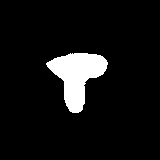}
\caption{Combined, $7\%$ noise}
\end{subfigure}
\begin{subfigure}{0.3\textwidth}
\centering
\includegraphics[width=0.85\textwidth]{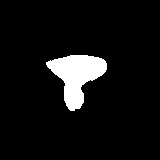}
\caption{Combined, $10\%$ noise}
\end{subfigure}
\begin{subfigure}{0.3\textwidth}
\centering
\includegraphics[width=0.85\textwidth]{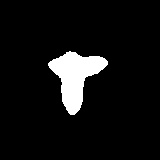}
\caption{Combined, $15\%$ noise}
\end{subfigure}
\caption{Reconstruction of a T-shaped object with $k = 6$, $\theta = 220^o$, and different noise levels added to the scattering data: $5\%$, $7\%$, $10\%$ and $15\%$.}
\label{T}
\end{figure}

\begin{figure}[tb]
\centering
\begin{subfigure}{0.3\textwidth}
\centering
\includegraphics[width=0.85\textwidth]{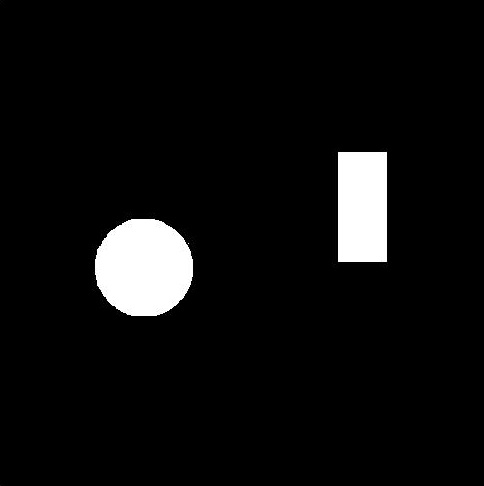}
\caption{True object}
\end{subfigure}
\begin{subfigure}{0.3\textwidth}
\centering
\includegraphics[width=0.85\textwidth]{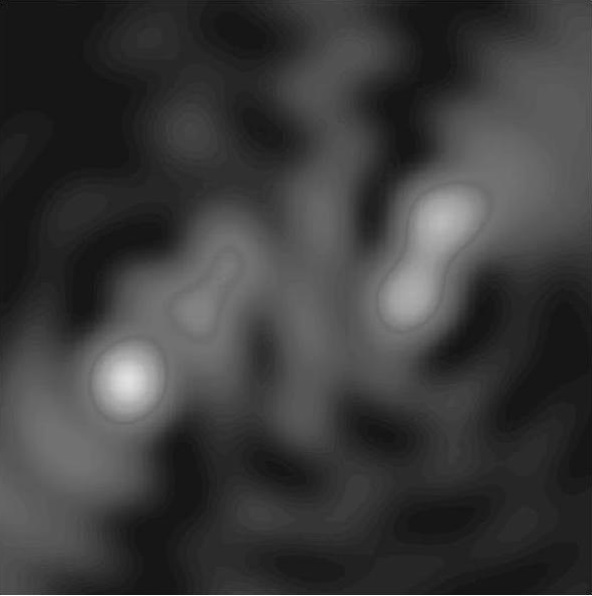}
\caption{Sampling, $5\%$ noise}
\end{subfigure}
\begin{subfigure}{0.3\textwidth}
\centering
\includegraphics[width=0.85\textwidth]{predict_7}
\caption{Combined, $5\%$ noise}
\end{subfigure}

\vspace{2mm}

\begin{subfigure}{0.3\textwidth}
\centering
\includegraphics[width=0.85\textwidth]{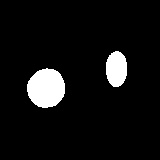}
\caption{Combined, $7\%$ noise}
\end{subfigure}
\begin{subfigure}{0.3\textwidth}
\centering
\includegraphics[width=0.85\textwidth]{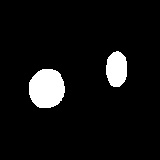}
\caption{Combined, $10\%$ noise}
\end{subfigure}
\begin{subfigure}{0.3\textwidth}
\centering
\includegraphics[width=0.85\textwidth]{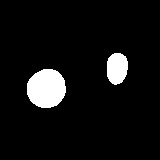}
\caption{Combined, $15\%$ noise}
\end{subfigure}
\caption{Reconstruction of a rectangle and a disk with $k = 6$, $\theta = 220^o$, and different noise levels added to the scattering data: $5\%$, $7\%$, $10\%$ and $15\%$.}
\label{er}
\end{figure}

Next, in Figures~\ref{pea}--\ref{er}, we test the method with different non-elliptical shapes and noisy scattering data associated with an incident direction different from which was used in the training process. This is to further assert its flexibility. Recall that we trained the neural network with $\theta = 90^o$ and $\theta = 45^o$. We reconstructed a peanut-shaped object, a T-shaped object and an object consisting of a disk and a rectangle using incident direction $\theta = 220^o$, see Figure \ref{pea}--\ref{er}. Again, we were able to get satisfactory results regardless of the circumstances.

For scattering data from a single incident direction and a fixed wave number, it is almost impossible not only for sampling methods but also for many other inversion methods  to reconstruct the shape and location of  extended objects. Therefore, these results show that the deep neural network has significantly improved the reconstruction capability of the imaging functional under lack of data. Moreover, they also show that the combined method is flexible, meaning it can work  well in some situations that are unseen from  the training process.

\section{Numerical study for experimental  data}
The results from simulated scattering data show that the proposed deep learning approach for OSM is capable of reconstructing different types of geometries in different setups. Therefore in this section, we demonstrate the effectiveness of our method when it comes to experimental data. 

We used the data sets for two-dimensional homogeneous objects provided by  Fresnel Institute~\cite{Belke2001}. Unlike the simulated data that are computed on an entire circle of radius 100 
these experimental data are measured on the arc of a circular sector  of angle $240^o$ (from $60^o$ to $300^o$ with step size $5^o$). The radius of the circle is 19 (about 760 mm).  Two  data sets are investigated: the first one named \textit{dielTM\_ dec4f.exp} deals with a de-centered circular cross section of radius $15$ mm, and the second one is \textit{rectTM\_cent.exp} concerning a centered rectangular cross section of dimensions $25.4 \times 12.7$ mm$^2$. We refer to~\cite{Belke2001} for detailed descriptions of the experimental setup.  

\begin{figure}[tb!]
\centering
\begin{subfigure}{0.3\textwidth}
\centering
\includegraphics[width=0.85\textwidth]{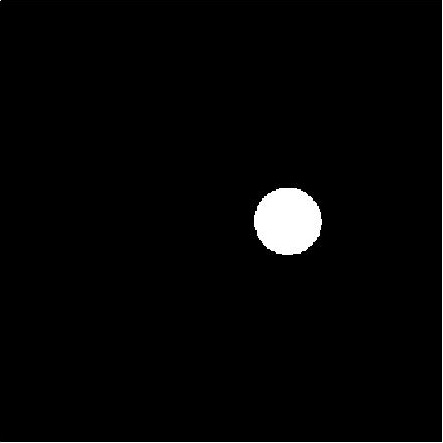}
\caption{True object}
\end{subfigure}
\begin{subfigure}{0.3\textwidth}
\centering
\includegraphics[width=0.85\textwidth]{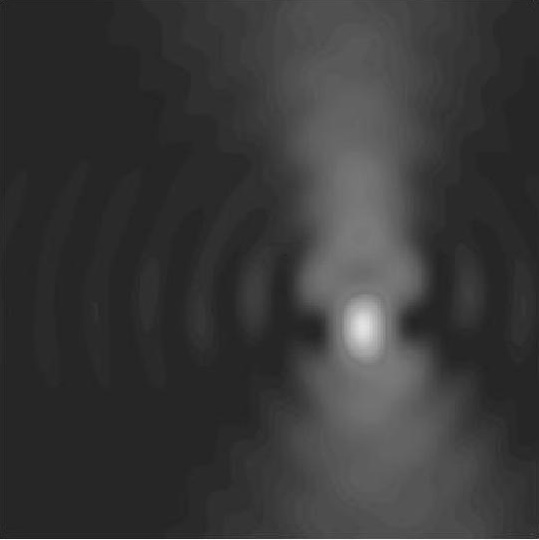}
\caption{Sampling method}
\end{subfigure}
\begin{subfigure}{0.3\textwidth}
\centering
\includegraphics[width=0.85\textwidth]{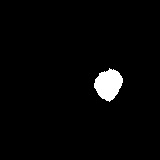}
\caption{Combined method}
\end{subfigure}
\caption{Reconstruction of a  circular object from experimental data at 8 GHz, $\theta = 90^o$.}
\label{real_diel}
\end{figure}

\begin{figure}[tb!]
\centering
\begin{subfigure}{0.3\textwidth}
\centering
\includegraphics[width=0.85\textwidth]{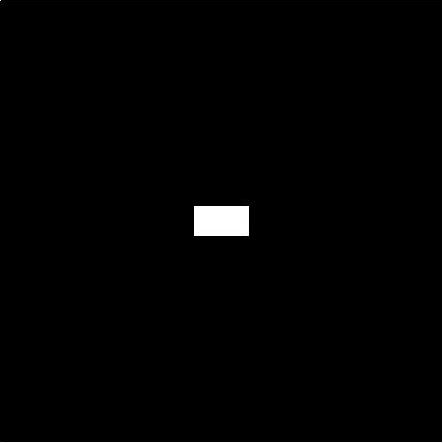}
\caption{True object}
\end{subfigure}
\begin{subfigure}{0.3\textwidth}
\centering
\includegraphics[width=0.85\textwidth]{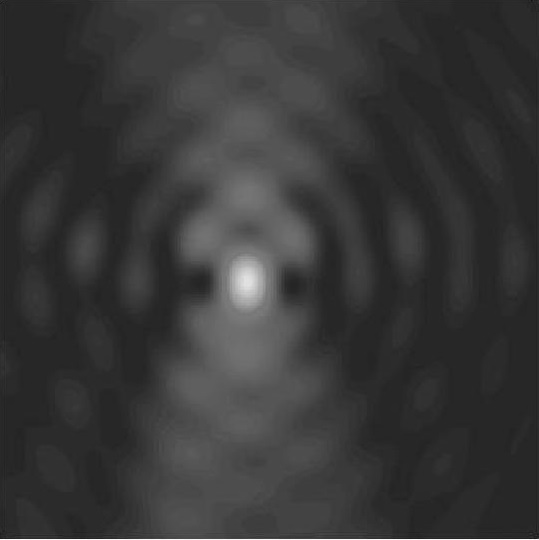}
\caption{Sampling method}
\end{subfigure}
\begin{subfigure}{0.3\textwidth}
\centering
\includegraphics[width=0.85\textwidth]{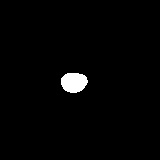}
\caption{Combined method}
\end{subfigure}
\caption{Reconstruction of a  rectangular object from experimental data at 8 GHz, $\theta = 90^o$.}
\label{real_rec}
\end{figure}

We rescale 40 millimeters to be 1 unit of length in our MATLAB simulations. No other adaptation  to  the  experimental  system  was necessary. The two scatterers were both tested at wave frequency 8 GHz (wave number $k$ is about 6.7) and with incident direction $\theta = 90^o$. See reconstructions using our combined method in Figure \ref{real_diel} and Figure \ref{real_rec}.

We see that the proposed method can handle  limited-aperture experimental data very well without any additional transfer learning. For the first object, we are able to achieve a reasonable reconstruction of both shape and location of the target. The second object is outside the training process, but the inverted result still can capture highly accurate information for its location and provide a pretty good estimate for its shape.

\bibliographystyle{amsplain}
\bibliography{ip-biblio}

%
%
%
%

\end{document}